\documentclass{amsart}
\usepackage{amssymb, amsthm, amscd, amsmath}

\usepackage{mathptm}
\usepackage{algorithmic}

\newtheorem{lemma}{Lemma}[section]
\newtheorem{corollary}[lemma]{Corollary}
\newtheorem{theorem}[lemma]{Theorem}
\newtheorem{proposition}[lemma]{Proposition}

\newtheorem{remark}[lemma]{Remark}

\begin{document}

\title{Rank one Maximal Cohen-Macaulay modules over\\ singularities of type $Y_1^3+Y_2^3+Y_3^3+Y_4^3$}
\author{Viviana \  \ \ {Ene}\ \ \ }
\address{Faculty of Mathematics and Informatics, Ovidius
University, Bd. Mamaia 124, 8700 Constanta, Romania}
\email{vivian@univ-ovidius.ro}
\author{\ \ \ Dorin\ \ \  {Popescu\ }}
\address{Institute of Mathematics, Bucharest University, P.O.
Box 1-764, \\ 70700 Bucharest, Romania}
\email{dorin@stoilow.imar.ro}


\begin{abstract}
We describe, by matrix factorizations, the rank one graded maximal
Cohen-Macaulay modules over the hypersurface
$Y_1^3+Y_2^3+Y_3^3+Y_4^3.$
\end{abstract}
\keywords{hypersurface ring, maximal Cohen-Macaulay modules,
Ulrich modules}
\subjclass[2000]{13C14,13H10, 13P10, 14J60}

\maketitle

\section{Introduction}

Let $R$ be a hypersurface ring, that is $R=S/(f)$ for a regular
local ring $(S,m)$ and $0\not = f\in m$. After Eisenbud
\cite{Ei2}, any maximal Cohen-Macaulay module  has a minimal free
resolution of periodicity 2 which is completely given by a matrix
factorization $(\phi,\psi)$, $\phi,\psi$ being square matrices
over $S$ such that $\phi \psi=\psi \phi=f I_n,$ for a certain
positive integer $n$. So, in order to describe the maximal
Cohen-Macaulay $R$-modules, it is enough to describe their matrix
factorizations (this we did for instance in \cite{EP} in order to
describe the maximal Cohen-Macaulay modules over singularities of
type $X^t+Y^3$). A different approach was used by Cipu, Herzog and
Popescu in \cite{CHP} to describe generalized Cohen-Macaulay
modules (see also \cite{C1} or \cite{CF}). A powerful method seems
to be also the lifting theory in the sense of
Auslander-Ding-Solberg \cite{ADS}, which was used in \cite{PP} in
order to complete Kn\"orrer Periodicity Theorem \cite{K} in char
$p>0$ (see also \cite{P}, \cite{C2}).\\

Let $R_n:=K[Y_1,\ldots,Y_n]/(f_n),$ where
$f_n=Y_1^3+Y_2^3+\ldots+Y_n^3$ and $K$ is an algebraically closed
field of characteristic $0.$ Using the classification of vector
bundles over elliptic curves obtained by Atiyah \cite{A}, C. Kahn
gives a "geometrically" description of the graded maximal
Cohen-Macaulay (briefly MCM) modules over $R_3$ and also describe
the Auslander-Reiten quivers of MCM over $R_3$ \cite{Ka}. His
method does not give the matrix factorizations of the
indecomposable MCM $R_3$-modules. In a recent paper \cite{LPP},
Laza, Pfister and Popescu use Atiyah classification to describe
the matrix factorizations of the graded, indecomposable, reflexive
modules over $R_3.$ They give canonical normal forms for the
matrix factorizations of these modules of rank one and show how
one may obtain the modules of rank $\geq 2$ using SINGULAR. Since
over the completion $K[[Y_1,Y_2,Y_3]]/(f_3)$ of $R_3,$ every
reflexive module is gradable (see \cite{Y1}), the authors obtain a
description of MCM-modules over $K[[Y_1,Y_2,Y_3]]/(f_3)$. \\ Now
we consider $n=4.$ In this case we do not have the support of
Atiyah classification used in the previous one, but we may give
the matrix factorizations for the rank one indecomposable MCM
modules over $R_4.$\\ Let $M$ be a MCM module over $R_4$ and let
$\mu(M)$ be the minimal number of generators of $M.$ By Corollary
1.3 of \cite{HK}, we obtain that $\mu(M)\in \{2,3\}.$ We shall
prove that there exists a finite number of indecomposable MCM
modules of rank one over $R_4.$ We note that, by \cite{LPP}, there
exists infinitely many indecomposable MCM modules of rank one over
$R_3.$ In \cite{HK}, Bruns showed that if $M$ is a MCM module over
a hypersurface ring, then rank $ M\geq (\dim R-1)/2.$ This implies
that there are no rank one MCM modules over $R_n,$ for $n\geq 5.$

\section{Rank one MCM modules over $R_4$ with two generators}

For every $a,b\in K$ with $a^3=b^3=-1$ and for every permutation
$(i\ j\ s)$ of the set $\{2,3,4\}$ with $i<j$ we denote:
 $$\varphi_{ij}(a,b)=
  \left(%
\begin{array}{cc}
    Y_1-aY_s & -(Y_i^2+bY_iY_j+b^2Y_j^2) \\
    Y_i-bY_j & Y_1^2+aY_1Y_s+a^2Y_s^2
  \end{array}%
\right)$$
 and

$$\psi_{ij}(a,b)=
  \left(%
\begin{array}{cc}
    Y_1^2+aY_1Y_s+a^2Y_s^2& (Y_i^2+bY_iY_j+b^2Y_j^2) \\
    -(Y_i-bY_j) & Y_1-aY_s
   \end{array}%
\right)$$

\begin{theorem}

$(\varphi_{ij}(a,b),\psi_{ij}(a,b))$ is a matrix factorization for
all $a,b\in K$ with $a^3=b^3=-1$ and $i,j\in \{2,3,4\}$ with
$i<j.$ The sets of graded MCM modules
 $$ \mathcal{M}=\{Coker\  \varphi_{ij}(a,b)|a,b,i,j\}$$

 and
 $$ \mathcal{N}=\{Coker\  \psi_{ij}(a,b)|a,b,i,j\} $$

 have the following properties:
\begin{description}
  \item[(i)] Every two generated, non free, graded MCM module is
  isomorphic with one of the modules of $\mathcal{M}\cup\mathcal{N}.$
  \item[(ii)] Every two different graded MCM modules from
  $\mathcal{M}\cup\mathcal{N}$ are not isomorphic.
  \item[(iii)] The modules of $\mathcal{M}$ are the syzygies and
  also the duals of the modules from  $\mathcal{N}.$
  \item[(iv)] The modules of $\mathcal{M}\cup\mathcal{N}$ are all
  of rank one.
\end{description}
\end{theorem}

\begin{proof}
 $(i)$ Obviously $(\varphi_{ij}(a,b),\psi_{ij}(a,b))$ is
a matrix factorization. Now let $(\varphi, \psi)$ be a reduced
$2\times 2$-matrix factorization of $f_4$ over
$K[Y_1,Y_2,Y_3,Y_4]$ with homogeneous entries. Then $\det \varphi
\det \psi=f_4^2$ and, since $f_4$ is irreducible, we have $\det
\varphi=\det \psi=f_4,$ after multiplication of a row of $\varphi$
and $\psi$ with some elements from $K^{\ast}.$ The matrix $\psi$
is the adjoint of $\varphi$, so it suffices to find $\varphi$ such
that $\det \varphi=f_4.$  After elementary transformations we may
suppose that the entries of the first column of $\varphi$ are
linear forms which must be linear independent since $f_4$ is
irreducible. So, applying some elementary transformations on the
matrix $\varphi,$ we may suppose that the entries of the first
column of $\varphi$ are of the form:
 $$\varphi_{11}=Y_1-a_{i_1}Y_{i_1}-a_{i_2}Y_{i_2}$$
 and
 $$\varphi_{21}=Y_i-b_{i_1}Y_{i_1}-b_{i_2}Y_{i_2}$$
 for some $a_{i_1},a_{i_2},b_{i_1},b_{i_2}\in K$,
$\{i,i_1,i_2\}=\{2,3,4\}$ and that the second column of $\varphi$
has the entries homogeneous forms of degree $2.$ Since $\det
\varphi=f_4$ we have that
$$f(a_{i_1}Y_{i_1}+a_{i_2}Y_{i_2},b_{i_1}Y_{i_1}+b_{i_2}Y_{i_2},Y_{i_1},Y_{i_2})=0.$$
This implies that $a_{i_1},a_{i_2},b_{i_1},b_{i_2}$ satisfy the
following identities:\begin{center}
\begin{tabular}{cc}
  $a_{i_1}^3+b_{i_1}^3+1$ &$=0,$ \\
  $ a_{i_2}^3+b_{i_2}^3+1$ &$ =0,$ \\
  $a_{i_1}^2a_{i_2}+b_{i_1}^2b_{i_2}$& $=0,$ \\
  $a_{i_1}a_{i_2}^2+b_{i_1}b_{i_2}^2$&$=0.$ \\
\end{tabular}
\end{center}
 If $b_{i_1}b_{i_2}\neq 0,$ then $ a_{i_1}a_{i_2}\neq 0$
and, from the last two equations, it results
$$(\frac{a_{i_1}}{b_{i_1}})^2=- \frac{b_{i_2}}{a_{i_2}}$$ and
$$(\frac{a_{i_2}}{b_{i_2}})^2=- \frac{b_{i_1}}{a_{i_1}},$$ so
$$(\frac{a_{i_1}}{b_{i_1}})^3= -1,$$ which contradicts the first
identity. Thus $b_{i_1}b_{i_2}= 0$ and $ a_{i_1}a_{i_2}= 0$. We
may suppose $b_{i_1}=0.$ It results: $$a_{i_1}^3=-1,\ a_{i_2}=0,\
b_{i_2}^3=-1.$$ We have obtained that $$\varphi_{11}=Y_1-aY_s$$
and $$\varphi_{21}=Y_i-bY_j,$$ where $a,b\in K,$ $a^3=b^3=-1,$ and
$ (i\ j\ s)$ is a permutation of the set $\{2,3,4\}$. It is clear
that we may transform the matrix such that $i<j.$ Let
$$\varphi=  \left(%
\begin{array}{cc}
  Y_1-aY_s & \gamma' \\
  Y_i-bY_j & \delta'
  \end{array}%
\right)$$
 where $\gamma', \delta'$ are homogeneous forms of degree
$2.$ Then we obtain that $\varphi$ and $\varphi_{ij}(a,b)$ define
the same MCM module as in (\cite{LOP}, Prop. 1.1).\\

 (ii) It is clear that no module of $\mathcal{M}$ is
isomorphic with one of $\mathcal{N}.$ The first Fitting ideal of
$\varphi_{ij}(a,b)$ is
Fitt$_1(\varphi_{ij}(a,b))=(Y_1-aY_s,Y_i-bY_j,Y_s^2,Y_j^2)$.
Suppose that $\varphi_{ij}(a,b)$ and $\varphi_{uv}(a',b')$ define
the same MCM module of $\mathcal{M}$. Then
$$\rm{Fitt}_1(\varphi_{ij}(a,b))=\rm{Fitt}_1(\varphi_{uv}(a',b'))$$ which
implies $$\varphi_{ij}(a,b)=\varphi_{uv}(a',b'),$$ as we can easy
check. Since the modules of $\mathcal{N}$ are the syzygies of
those of $\mathcal{M}$ it results that any two different modules
of $\mathcal{N}$ are not isomorphic. (3) and (4) follows as in
(\cite{LPP}, Theorem 3.1).
\end{proof}

\begin{remark}{\em
We note that every matrix factorization of a two generated, non
free, graded MCM module over $R_4$ is the tensor product of the
matrix factorizations of $Y_1^3+Y_s^3$ and $Y_i^3+Y_j^3$ (see
\cite{Y2}).}
\end{remark}

\section{Rank one MCM modules over $R_4$ with three generators}

 Let $M$ be a rank one MCM module over $R_4$ with three
generators and let $(\varphi,\psi)$ be a matrix factorization of
$M.$ We may suppose $\det \varphi=f_4$ (if necessary replacing $M$
by its first syzygy). Thus the entries of $\varphi$ are linear
forms.

\begin{lemma}
Let $\alpha,\beta,\gamma,\delta$ be independent linear forms in
$K[Y_1,Y_2,Y_3,Y_4]$ such that $f_4\in
(\alpha,\beta)\cap(\gamma,\delta).$ Then there exists some linear
forms $m,n,w,t$ such that $$\det\left(%
\begin{array}{ccc}
  0 & \alpha & \beta \\
\gamma & m & n \\
 \delta & w & t
\end{array}%
\right)=f_4.$$
\end{lemma}

\begin{proof}
Since $f_4\in (\alpha,\beta)$ there exist non unique  $2-$forms
$\eta_1,\eta_2$ such that \begin{equation} f_4=\alpha \eta_1+\beta
\eta_2.\end{equation} $ \eta_1,\eta_2$ can be expressed as:
\begin{equation}\eta_1=\eta_{11}\alpha +\eta_{12}\beta
+\eta_{13}\gamma+\eta_{14}\delta \end{equation}
\begin{equation}\eta_2=\eta_{21}\alpha+\eta_{22}\beta+\eta_{23}\gamma+\eta_{24}\delta,\end{equation}
where $\eta_{ij}$ are linear forms, since
$\alpha,\beta,\gamma,\delta$ are independent and so generate the
linear form space. By hypothesis $f_4\in (\gamma,\delta)$ so
$$\alpha \eta_1+\beta \eta_2\equiv 0 \pmod {(\gamma,\delta)},$$
which implies $\alpha \eta_1 \equiv 0 \pmod
{(\beta,\gamma,\delta)}.$ But $\alpha$ is not contained in the
prime ideal $(\beta,\gamma,\delta)$. It results that
$$\eta_1\equiv 0 \pmod {(\beta,\gamma,\delta)},$$ thus we may take
$\eta_{11}=0. $ Replacing the expressions of $\eta_1$ and $\eta_2$
in the equality (1) we get $$(\eta_{12}+\eta_{21})\alpha \beta
+\eta_{22}\beta^2 \in(\gamma,\delta).$$ Since $\beta \notin
(\gamma,\delta)$ we deduce that
\begin{equation}(\eta_{12}+\eta_{21})\alpha +\eta_{22}\beta
\in(\gamma,\delta).\end{equation} This implies that
$$\eta_{22}\beta \equiv 0 \pmod {(\alpha,\gamma,\delta)}.$$
Moreover, we have $\eta_{22}\equiv 0 \pmod
{(\alpha,\gamma,\delta)}$. It follows that there exists
$\lambda_1,\lambda_2,\lambda_3\in K$ such that
$$\eta_{22}=\lambda_1\alpha+\lambda_2\gamma+\lambda_3\delta.$$ By
the relation (4) we have that
$$\eta_{12}+\eta_{21}+\lambda_1\beta\equiv 0
\pmod{(\gamma,\delta)}$$ so $$\eta_{21}\equiv
-\eta_{12}-\lambda_1\beta \pmod{(\gamma,\delta)}.$$ Therefore we
may write $\eta_2$ in the following form:
$$\eta_2=-\eta_{12}\alpha +\eta_{23}'\gamma +\eta_{24}'\delta.$$
Denote $\eta_1'=\eta_1-\eta_{12}\beta$ and
$\eta_2'=\eta_2+\eta_{12}\alpha.$ Then
$$f_4=\alpha\eta_1'+\beta\eta_2'$$ and $$\eta_1',\eta_2'\equiv 0
\pmod {(\gamma,\delta)}.$$ Thus we may find some linear forms with
the required property.
\end{proof}

 For $1\leq i\leq 4$, let $\mathcal{L}_i$ be the set of
the linear forms $Y_i-aY_j,$ where $a\in K,\ a^3=-1$ and $j\in
\{1,2,3,4\},\ j>i.$

\begin{proposition}\label{main}
Let $M$ be a three generated, rank one, graded MCM mo\-dule over
$R_4$. Then there exist some independent linear forms
$\alpha,\beta,\gamma,\delta$ with $\alpha,\gamma \in
\mathcal{L}_1,\ \beta \in \mathcal{L}_j$ and $\delta \in
\mathcal{L}_i$ for some $i,j\geq 2$ and there exist $m,n,w,t$
linear forms such that $$\varphi=\left(%
\begin{array}{ccc}0 & \alpha &
\beta
\\ \gamma & m & n\\
 \delta &w & t
\end{array}%
\right)$$
 and its adjoint matrix, $\psi,$ form a matrix
factorization of $M.$
\end{proposition}
\begin{proof}
As rank $M=1$, every matrix factorization $(\varphi,\psi)$ of $M$
has $\det \varphi=f_4.$ Since $f_4\in (Y_1+Y_2,Y_3+Y_4)$, we
obtain that $\varphi$ has a generalized zero (see \cite{Ei1}). By
elementary transformations $\varphi$ can be arranged in
the form $$\varphi=\left(%
\begin{array}{ccc} 0 & \alpha & \beta
\\ \gamma & m & n\\
 \delta &w & t
\end{array}%
\right).$$
 As in the two generated case we obtain
$$\alpha=Y_1-aY_{j_1},\beta = Y_j-bY_{j_2}, \gamma = Y_1-cY_{i_1},
\delta= Y_i-dY_{i_2},$$ where $(j,j_1,j_2)$ and $(i,i_1,i_2)$ are
permutations of the set $\{2,3,4\}$ such that $j<j_2$ and $i<i_2$,
that is $\alpha,\gamma \in \mathcal{L}_1,\ \beta \in
\mathcal{L}_j$ and $\delta \in \mathcal{L}_i.$ We shall prove that
since $\det \varphi=f_4$ we must have $\alpha,\beta,\gamma,\delta$
linear independent. We have the following possibilities to choose
$\varphi:$
\begin{description}
  \item[(i)]$A=\left(%
\begin{array}{ccc}
   0 & Y_1-aY_4 & Y_2-bY_3 \\
 Y_1-cY_2 & \star & \star \\
    Y_3-dY_4 & \star & \star \
 \end{array}%
\right)$
  \item[(ii)] $A^t$
  \item[(iii)]$B=\left(%
\begin{array}{ccc}
   0 & Y_1-aY_3 & Y_2-bY_4 \\
 Y_1-cY_2 & \star & \star \\
    Y_3-dY_4 & \star & \star \
\end{array}%
\right)$
 \item[(iv)] $B^t$
 \item[(v)] $C=\left(%
\begin{array}{ccc}
   0 & Y_1-aY_4 & Y_2-bY_3 \\
 Y_1-cY_3 & \star & \star \\
    Y_2-dY_4 & \star & \star \
\end{array}%
\right)$
  \item[(vi)] $C^t$
  \item[(vii)] $D=\left(%
\begin{array}{ccc}
   0 & Y_1-aY_2 & Y_3-bY_4 \\
 Y_1-cY_2 & \star & \star \\
    Y_3-dY_4 & \star & \star \
  \end{array}%
\right)$
  \item[(viii)] $E=\left(%
\begin{array}{ccc}
   0 & Y_1-aY_3 & Y_2-bY_4 \\
 Y_1-cY_3 & \star & \star \\
    Y_2-dY_4 & \star & \star \
 \end{array}%
\right)$
  \item[(ix)] $F=\left(%
\begin{array}{ccc}
   0 & Y_1-aY_4 & Y_2-bY_3 \\
 Y_1-cY_4 & \star & \star \\
    Y_2-dY_3 & \star & \star \
\end{array}%
\right)$
\end{description}

 We shall give the proof only for the first case. The
others are similar. \\  Let $ \varphi=\left(%
\begin{array}{ccc}
   0 & Y_1-aY_4 & Y_2-bY_3 \\
 Y_1-cY_2 & m & n \\
    Y_3-dY_4 & w & t \
 \end{array}%
\right).$ Since $\det \varphi =f_4$ we obtain:
  $$-(Y_1-aY_4)((Y_1-cY_2)t-(Y_3-dY_4)n)+(Y_2-bY_3)((Y_1-cY_2)w-(Y_3-dY_4)m)=$$
  $$(Y_1-aY_4)(Y_1^2+aY_1Y_4+a^2Y_4^2)+(Y_2-bY_3)(Y_2^2+bY_2Y_3+b^2Y_3^2).$$
 This equality is equivalent with
$$(Y_1-aY_4)((Y_1^2+aY_1Y_4+a^2Y_4^2)+ (Y_1-cY_2)t-(Y_3-dY_4)n)=$$
$$(Y_2-bY_3)((Y_1-cY_2)w-(Y_3-dY_4)m-(Y_2^2+bY_2Y_3+b^2Y_3^2)).$$
 It results that there exists a linear form $\eta$ such
that $$ (Y_1^2+aY_1Y_4+a^2Y_4^2)+
(Y_1-cY_2)t-(Y_3-dY_4)n=\eta(Y_2-bY_3).$$ Put
$$Y_1=bcd,Y_2=bd,Y_3=d\ \rm{and}\  Y_4=1$$ in the above equality. It
follows $$(bcd)^2+abcd+a^2=0,$$ which gives $a\neq bcd.$ This
condition means exactly that $\alpha,\beta,\gamma,\delta$ are
linearly independent. Thus $bcd=\epsilon a$,where $\epsilon$ is in
$K,$
 $\epsilon^3=1$ and $\epsilon\neq 1.$ An example of such $A$
 is given by:
$$A=\left(%
\begin{array}{ccc}
   0 & Y_1-aY_4 & Y_2-bY_3 \\
 Y_1-cY_2 & -b^2Y_3-abc^2\epsilon^2Y_4 & b^2c^2Y_3-abc\epsilon^2Y_4 \\
    Y_3-dY_4 & c^2Y_2+bc^2Y_3+acY_4 & -Y_1-cY_2-aY_4 \
 \end{array}%
\right).$$ Then $A$ and its adjoint, $A^{\ast},$ form a matrix
factorization of $f_4$ (by the following lemma we see that always
$A$ can be supposed
of the above form after some elementary transformations).\\
The condition of linear independence of
$\alpha,\beta,\gamma,\delta,$ in the case {\bf (iii)} is $ad\neq
bc,$ that is $ad=\epsilon bc.$ Then $$B=\left(%
\begin{array}{ccc}
   0 & Y_1-aY_3 & Y_2-bY_4 \\
 Y_1-cY_2 & a^2cY_3+(abc^2+a^2cd)Y_4 & a^2Y_3-a^2d\epsilon Y_4 \\
    Y_3-dY_4 & c^2Y_2+acY_3+bc^2Y_4 & -Y_1-cY_2-aY_3 \
\end{array}%
\right)$$ and its adjoint matrix, $B^{\ast},$ form a matrix
factorization of $f_4.$\\
The condition of linear independence of
$\alpha,\beta,\gamma,\delta,$ in the case {\bf (v)} is $ab\neq
cd,$ that is $ab=\epsilon cd.$ Then $$C=\left(%
\begin{array}{ccc}
   0 & Y_1-aY_4 & Y_2-bY_3 \\
 Y_1-cY_3 & -Y_2-bY_3-dY_4 & -b^2c^2Y_3+bc^2d\epsilon^2 Y_4\\
    Y_2-dY_4 & b^2c^2Y_3-bc^2d\epsilon^2 Y_4  & -Y_1-cY_3-aY_4 \
\end{array}%
\right)$$ and its adjoint matrix, $C^{\ast},$ form a matrix
factorization of $f_4.$\\
For the last three cases we obtain that
$\alpha,\beta,\gamma,\delta$ are linear independent if and only if
$a\neq c$ and $b\neq d.$ Then the pairs $(D,D^{\ast}),\
(E,E^{\ast})$ and $(F,F^{\ast})$ are matrix factorizations, where
$$D=\left(%
\begin{array}{ccc}
   0 & Y_1-aY_2 & Y_3-bY_4 \\
 Y_1-cY_2 & -Y_3-(b+d)Y_4& 0 \\
    Y_3-dY_4 & 0 & -Y_1-(a+c)Y_2 \
  \end{array}%
\right),$$ $$E=\left(%
\begin{array}{ccc}
   0 & Y_1-aY_3 & Y_2-bY_4 \\
 Y_1-cY_3 & -Y_2-(b+d)Y_4 & 0 \\
    Y_2-dY_4 & 0 &-Y_1-(a+c)Y_3 \
 \end{array}%
\right)$$ and $$F=\left(%
\begin{array}{ccc}
   0 & Y_1-aY_4 & Y_2-bY_3 \\
 Y_1-cY_4 & -Y_2-(b+d)Y_3 & 0 \\
    Y_2-dY_3 & 0 & -Y_1-(a+c)Y_4 \
\end{array}%
\right).$$
\end{proof}

The next Lemma will show that every three generated, rank one,
non-free graded MCM module over $R_4$ is isomorphic with a module
given by one of the above matrix factorizations.

\begin{lemma}
If $\alpha,\beta,\gamma,\delta$ are independent linear forms as in
the above Proposition and $$\varphi=\left(%
\begin{array}{ccc}
  0 & \alpha & \beta \\
\gamma & m & n \\
 \delta & w & t
\end{array}%
\right),\ \varphi'=\left(%
\begin{array}{ccc}
  0 & \alpha & \beta \\
\gamma & m' & n' \\
 \delta & w' & t'
\end{array}%
\right)$$
 then $\mbox{Coker\ } \ \varphi \cong \mbox{Coker\ }\varphi'.$
\end{lemma}

\begin{proof} Let $\eta$ and $\nu$ be two homogeneous forms of degree
2 such that $f_4=\alpha\eta+\beta\nu.$ It results that $$\alpha n
\delta + \beta\gamma w- \delta\beta m - \alpha\gamma t=
\alpha\eta+\beta\nu,$$ that is $$\alpha(n \delta - \gamma t-\eta)=
\beta(\nu - \gamma w +\delta m). $$ Therefore we obtain the
following equalities: \begin{equation}n \delta - \gamma t-\eta
=\theta \beta\end{equation}  and \begin{equation}\delta m - \gamma
w + \nu = \theta \alpha,\end{equation}  for some linear form
$\theta.$ In the same way we obtain that there exists a linear
form $\theta'$ such that
\begin{equation}
n' \delta - \gamma t'-\eta =\theta' \beta\end{equation} and
 \begin{equation}\delta m'- \gamma w' + \nu = \theta' \alpha,\end{equation}
 Subtracting the identities (5) and (7) we obtain:
 \begin{equation}(n-n')\delta-(t-t')\gamma=(\theta-\theta')\beta.\end{equation}
  Since $\beta \notin (\gamma,\delta)$
 it follows that there exist $a,b\in K$ such that $\theta-\theta'=a\delta +b\gamma.$
Replacing in the equation (9) we get:
$$(n-n'-a\beta)\delta=(t-t'+b\beta)\gamma.$$ Thus there exists
$c\in K$ such that
\begin {equation}
t'=t+b\beta-c\delta
\end{equation}
and
\begin{equation}
n'=n-a\beta -c\gamma.
\end{equation}
 Starting with the equations (6) and (8) we obtain analogously
that there exists $c'\in K$ such that
\begin{equation}
m'=m-a\alpha-c'\gamma
\end{equation}
and
\begin{equation}
w'=w+ b\alpha -c'\delta.
\end{equation}
The last four equalities show that $\varphi'$ is obtained from
$\varphi$ after some elementary transformations and so prove our
Lemma.
\end{proof}

 From now on, the most difficult task is to decide which of the
modules given by the matrix factorizations defined in the proof of
Proposition \ref{main} are isomorphic. We recall that two
matrices, $\varphi$ and $\varphi',$ define the same module over
$R_4$ (i.e. Coker $\varphi\simeq$ Coker $\varphi'$), if and only
if they are equivalent, that is there exist $U$ and $V$ two square
matrices with entries in $K[Y_1,\ldots,Y_n]$ such that
$\varphi'=U\varphi V$ and $\det(U)=\det(V)=1$ (see \cite{Ei2}). In
this case we denote $\varphi
\sim \varphi'.$\\
The proof of the main theorem of this section will be done with
the help of the computer algebra system SINGULAR \cite{GPS}.\\ For
$a,b,c,d,\epsilon\in K$ such that $a^3=b^3=c^3=d^3=-1,
\epsilon^3=1,\epsilon\neq 1$ and
$bcd=\epsilon a,$ we set $$\alpha(b,c,d,\epsilon)=\left(%
\begin{array}{ccc}
   0 & Y_1-aY_4 & Y_2-bY_3 \\
 Y_1-cY_2 & -b^2Y_3-abc^2\epsilon^2Y_4 & b^2c^2Y_3-abc\epsilon^2Y_4 \\
    Y_3-dY_4 & c^2Y_2+bc^2Y_3+acY_4 & -Y_1-cY_2-aY_4 \
 \end{array}%
\right)$$ and $$\beta(b,c,d,\epsilon)=\alpha(b,c,d,\epsilon)^t,$$
that is the transpose of $\alpha(a,b,c,d).$  We know from the
proof of the Proposition \ref{main} that $(\alpha(b,c,d,\epsilon),
\alpha(b,c,d,\epsilon)^{\ast})$ and
$(\beta(b,c,d,\epsilon),\beta(b,c,d,\epsilon)^{\ast})$ are matrix
factorizations of $f_4$.\\ For $a,b,c\in K,$ distinct
roots of $-1,$ and $\epsilon$ as above, we set $$\eta(a,b,c,\epsilon)=\left(%
\begin{array}{ccc}
   0 & Y_1+Y_2 & Y_3-aY_4 \\
 Y_1+\epsilon Y_2 & -Y_3+cY_4 & 0 \\
    Y_3-bY_4 & 0 &-Y_1-\epsilon^2Y_2 \
 \end{array}%
\right),$$ $$\vartheta(a,b,c)=\left(%
\begin{array}{ccc}
   0 & Y_1+Y_3 & Y_2-aY_4 \\
 Y_1-a^2bY_3 & -Y_2+cY_4 & 0 \\
    Y_2-bY_4 & 0 &-Y_1+ab^2Y_3 \
 \end{array}%
\right).$$ These matrices are of the type  $D$ and $E.$ Thus,
every matrix forms with its adjoint a matrix factorization of
$f_4$.

\begin{theorem}
 Let
$$\mathcal{M}=\{\mbox{Coker\ } \alpha(b,c,d,\epsilon),\ \mbox{Coker\ } \beta(b,c,d,\epsilon)\
|\ b,c,d,\epsilon\in K, \ $$ $$b^3=c^3=d^3=-1,\ bcd=\epsilon a,\
\epsilon^3=1,\epsilon\neq 1 \}$$ and $$\mathcal{N}=\{\mbox{Coker\
}\eta(a,b,c,\epsilon),\ \mbox{Coker\ }\vartheta(a,b,c),\
 |\ \epsilon^3=1,\ \epsilon\neq 1 $$ $$\mbox{\ and\ } (a,b,c) \mbox{\ is\ a\
permutation\ of\ the\ roots\ of\ }-1\}.$$ Then the sets
$\mathcal{M},\mathcal{N}$ of rank one, three generated, MCM graded
$R_4$-modules have the following properties:
\begin{description}
    \item[(i)] every three generated, rank one, non-free, graded
    MCM $R_4-$module is isomorphic with one module from $\mathcal{M}\cup \mathcal{N}.$
    \item[(ii)] if $M=$ Coker $\alpha(b,c,d,\epsilon)$ (or
    $M=$ Coker $\beta(b,c,d,\epsilon)$) belongs to $\mathcal{M}$ and $N\in
    \mathcal{M},$ then $N\simeq M$ if and only if $N=$ Coker $\alpha (b\epsilon,c
    \epsilon,d\epsilon,\epsilon^2)$ (or $N=\beta(b\epsilon,c
    \epsilon,d\epsilon,\epsilon^2)$).
    \item[(iii)] any two different modules from $\mathcal{N}$
    are not isomorphic.
    \item[(iv)] any module of $\mathcal{N}$ is not isomorphic with
    some module of $\mathcal{M}.$
\end{description}
\end{theorem}

\begin{proof}
For the beginning we shall prove that any module of the type
$B,B^t, C$ and $C^t$ of the proof of Proposition \ref{main} is
isomorphic with one of type $A$ or $A^t.$ This can be done using
SINGULAR. For instance, to establish that the modules of type $B$
are isomorphic with modules of type
$A^t,$ we use the following procedure (see \cite{LPP}, Lemma 5.1):\\

\noindent{\tt
LIB"matrix.lib"; \\
option(redSB);\\
proc isomorf(matrix X, matrix Y)\\
\{ \\ matrix U[3][3]=u(1..9);\\ matrix V[3][3]=v(1..9);\\ matrix
C=U*X-Y*V;\\ ideal I=flatten(C);\\ ideal
I1=transpose(coeffs(I,y(1)))[2];\\ ideal
I2=transpose(coeffs(I,y(2)))[2];\\ ideal
I3=transpose(coeffs(I,y(3)))[2];\\ ideal
I4=transpose(coeffs(I,y(4)))[2];\\ ideal
J=I1+I2+I3+I4+ideal(det(U)-1,det(V)-1);\\ ideal
L=std(J);\\
L;\\ \}}\\

\noindent We apply this procedure for the matrices $A^t$ and
$B:$\\

\noindent{\tt ring R=0,(u(1..9),v(1..9),y(1..4),x,a,b,c,d,m,n,p,q,y),lp;\\
ideal
F=a3+1,b3+1,c3+1,d3+1,x*a-b*c*d,x2+x+1,m3+1,n3+1,p3+1,q3+1,\\m*q-y*n*p,y2+y+1;\\
qring Q=std(F);\\ matrix
A[3][3]=0,y(1)-a*y(4),y(2)-b*y(3),y(1)-c*y(2),-b2*y(3)\\-a*b*c2*x2*y(4),
b2*c2*y(3)-a*b*c*x2*y(4),y(3)-d*y(4),c2*y(2)\\+b*c2*y(3)+a*c*y(4),-y(1)-c*y(2)-a*y(4);\\
matrix
B[3][3]=0,y(1)-m*y(3),y(2)-n*y(4),y(1)-p*y(2),m2*p*y(3)\\+m*n*p2*y(4)+m2*p*q*y(4),m2*y(3)-m2*q*y*y(4),y(3)-q*y(4),p2*y(2)\\+m*p*y(3)+n*p2*y(4),
-y(1)-p*y(2)-m*y(3);\\ // Now we test the equivalence between the
matrices transpose(A) // and B\\ isomorf(transpose(A),B); }\\

We obtain that $A^t$ and $B$ are equivalent if and only if
\begin{equation}\label{1atb}
d^2-dqy-dq+q^2y =0
\end{equation}
\begin{equation}\label{2atb}
c-dpq^2y-dpq^2=0
\end{equation}
\begin{equation}\label{3atb}
 b-dnq=0
 \end{equation}
 \begin{equation}\label{4atb}
a+dnpq^2+npy+np=0
\end{equation}
 and
 \begin{equation}\label{5atb}
ab^2c^2d^2+dq^2y-1=0.
\end{equation}

If $m,n,p,q$ and $y$ are fixed such that $m^3=n^3=p^3=q^3=-1,\
y^2+y+1=0$ and $mq=npy,$ then  we may obtain $a,b,c,d$ and $x$
such that the above equations are satisfied and
$a^3=b^3=c^3=d^3=-1,\ x^2+x+1=0,\ bcd=ax.$\\ For instance, we may
take $$a=-np,\ b=ynq^2,\ c=p {\rm \ and\ } d=qy.$$ In this case,
if $\lambda\in K$ verifies the equation $9\lambda^3=8mp^2,$ one
finds the matrices $$U=\left(%
\begin{array}{ccc}
  0 & -\frac{1}{2}p\lambda &  \frac{1}{2}p^2\lambda\\
  \frac{3}{2}m^2\lambda & 0 & 0 \\
  0 & -\frac{1}{2}p^2\lambda & \lambda \\
\end{array}\right){\rm \ and\ }V=\left(%
\begin{array}{ccc}
  0 & \frac{3}{2}m^2\lambda & 0\\
  -\frac{1}{2}p\lambda & 0 & -\frac{1}{2}p^2\lambda  \\
  \frac{1}{2}p^2\lambda  & 0 & \lambda \\
\end{array}\right),$$ such that $UA^t=BV.$\\ With the same
procedure, we obtain that every matrix of type $C$ wich depends on
$m,n,p,q$ and $y$ is equivalent with the transpose of a matrix of
type $A$ depending on $a,b,c,d$ and $x,$ where $$a=-n^2pqy^2,\
b=n,\ c=n^2p,\ d=n^2q\rm{\ and\ }x=y.$$
 Now we study the equivalence of the matrices of type $D,E,F.$ Let $(a,b,c)$ and
$(p,q,r)$ be two permutations of the third roots of $-1$ and $$D((a,b,c),(p,q,r))=\left(%
\begin{array}{ccc}
   0 & Y_1-aY_2 & Y_3-pY_4 \\
 Y_1-cY_2 & -Y_3+rY_4& 0 \\
    Y_3-qY_4 & 0 & -Y_1+bY_2 \
  \end{array}%
\right).$$ One can apply elementary transformations on the columns
and on the rows of $D((a,b,c),(p,q,r))$ to obtain
$$D((a,b,c),(p,q,r))\sim D((c,a,b),(q,r,p))\sim
D((b,c,a),(r,p,q)).$$ We deduce similar equivalences for the
matrices of type $E$ and $F.$ This means that we may restrict our
study to the matrices
$$\eta(p,q,r,\epsilon)=D((-1,-\epsilon,-\epsilon^2),(p,q,r)),$$ $$
E((-1,-\epsilon,-\epsilon^2),(p,q,r)),\
F((-1,-\epsilon,-\epsilon^2),(p,q,r)),$$ where $\epsilon$ is a
root $\neq 1$ of $1$ in $K.$ \\
Using the procedure {\tt isomorf(matrix X, matrix Y)}, we get that
any two different matrices $D((-1,-x,-x^2),(a,b,c))$ and
$D((-1,-y,-y^2),(p,q,r))$ are not equivalent.\\ Now, let us
consider the matrices $D((-1,-x,-x^2),(a,b,c))$ and
$E((-1,-y,-y^2),\\(p,q,r)).$ Applying our procedure, it results
that they are equivalent if and only if $$y+pq^2=0,\ b-q=0,\
a+p+q=0 {\rm\ and\ } x^2-pq^2+1=0.$$ Since $-1, pq^2,p^2q$ are the
solutions of the equation $x^3+1=0,$ we get $x^2=-1+pq^2=-p^2q.$
Thus, if $y=-pq^2,$ then $E((-1,-y,-y^2),(p,q,r))$ is equivalent
with $D((-1,pq^2,p^2q),(-p-q,q,p)).$ If $y=-p^2q$ (the only left
case!) then $E((-1,-y,-y^2),(p,q,r))$ is $\vartheta(p,q,r).$ Until
now we have obtained that the matrices which define the modules of
the set $\mathcal{N}$ are pairwise non-equivalent. Finally, we
find that the matrices of the form $F((-1,-y,-y^2),(p,q,r)),$
where $y=-pq^2,$ are equivalent with some matrices of type $D,$
and those $F((-1,-y,-y^2),(p,q,r)),$ where $y=-p^2q,$ are
equivalent with some matrices of the type
$\vartheta(a,b,c),$ with $b=q$ and $a+p+q=0.$\\
 So we have proved the parts {\bf(i)} and {\bf (iii)} of the theorem. For the rest, one can use the
procedure {\tt isomorf(matrix X, matrix Y)}, as in the previous
part. For instance, to prove {\bf (iv)}:\\

\noindent{\tt ring R=0,(u(1..9),v(1..9),y(1..4),x,a,b,c,d,w,p,q),lp;\\
ideal F=a3+1,b3+1,c3+1,d3+1,x*a-b*c*d,x2+x+1,p3+1,q3+1,w2-w+1;\\
qring Q=std(F);\\ matrix
A[3][3]=0,y(1)-a*y(4),y(2)-b*y(3),y(1)-c*y(2),-b2*y(3)\\-a*b*c2*x2*y(4),b2*c2*y(3)-a*b*c*x2*y(4),y(3)-d*y(4),c2*y(2)\\+b*c2*y(3)+a*c*y(4),-y(1)-c*y(2)-a*y(4);\\
matrix
D[3][3]=0,y(1)+y(2),y(3)-p*y(4),y(1)-w*y(2),-y(3)-p*y(4)\\-q*y(4),0,y(3)-q*y(4),0,-y(1)-w2*y(2);\\
// Now we test the equivalence between the
matrices A and D\\
isomorf(A,D);\\
// Now we test the equivalence between the
matrices transpose(A) // and D\\
isomorf(transpose(A),D);}\\

\noindent In both cases we obtain:\\

 \noindent {\tt L[1]=1}\\

\noindent This proves that there is no module of type $D$ which is
isomorphic with a module of the set $\mathcal{M}.$ Analogously we
may check that there is no module of type $E$ which
is isomorphic with a module of the set $\mathcal{M}.$ This shows {\bf (iv)}.\\
Finally, for the part {\bf (ii)}, we apply the procedure {\tt
isomorf(matrix X, matrix Y)} for the matrices $\alpha(b,c,d,x)$
and $\alpha(n,p,q,y).$ We obtain that these two matrices are
equivalent if and only if the following equations are satisfied:

\begin{equation}\label{1aa}
d^2-dqy-dq+q^2y =0
\end{equation}
\begin{equation}\label{2aa}
c+dpq^2=0
\end{equation}
\begin{equation}\label{3aa}
b+dnq^2=0
\end{equation}
\begin{equation}\label{4aa}
a+dnpy+npq=0
\end{equation}
\begin{equation}\label{5aa}
ab^2c^2d^2+dq^2y-1=0.
\end{equation}
 From the equation (\ref{1aa}) we obtain: (i) $\frac{d}{qy^2}=y$ or (ii) $\frac{d}{qy^2}=y^2.$
In the case (i), it follows $d=q$ and, from the above equations,
we obtain $c=p,\ b=n$ and $a=y^2npq=my^3=m.$ The equation
(\ref{5aa}) is obviously verified. Thus in the first case we get
that $\alpha(b,c,d,x)=\alpha(n,p,q,y).$ In the second one, using
the equations (\ref{1aa})-(\ref{5aa}), we obtain $$d=qy,\ c=py,\
b=ny{\rm\ and\ }x=y^2.$$ To finish the proof of {\bf(ii)} we apply
the procedure {\tt isomorf(matrix X, matrix Y)} for the
matrices $\beta(b,c,d,x)$ and $\alpha(n,p,q,y):$ \\

{\tt \noindent ring R=0,(u(1..9),v(1..9),y(1..4),x,a,b,c,d,m,n,p,q,y),lp;\\
ideal
F=a3+1,b3+1,c3+1,d3+1,x*a-b*c*d,x2+x+1,m3+1,n3+1,p3+1,\\q3+1,y*m-n*p*q,y2+y+1;\\
qring Q=std(F);\\ matrix
A[3][3]=0,y(1)-a*y(4),y(2)-b*y(3),y(1)-c*y(2),-b2*y(3)\\-a*b*c2*x2*y(4),b2*c2*y(3)-a*b*c*x2*y(4),y(3)-d*y(4),c2*y(2)\\+b*c2*y(3)+a*c*y(4),-y(1)-c*y(2)-a*y(4);\\
matrix
AA[3][3]=0,y(1)-m*y(4),y(2)-n*y(3),y(1)-p*y(2),-n2*y(3)\\-m*n*p2*y2*y(4),n2*p2*y(3)-m*n*p*y2*y(4),y(3)-q*y(4),p2*y(2)\\+n*p2*y(3)+m*p*y(4),-y(1)-p*y(2)-m*y(4);\\
isomorf(transpose(A),AA); }\\

\noindent We obtain:\\

 \noindent {\tt L[1]=1}\\

\noindent This shows that no matrix of type $A$ is equivalent with
one of type $A^t.$
\end{proof}

The three generated, rank one, MCM modules over $R_4$ are linear
MCM or Ulrich modules (see \cite{BH}, \cite{BHS}). Thus, from the
above theorem we obtain:
\begin{corollary}\label{victor}
There are $72$ isomorphism classes of Ulrich modules of rank one
over the ring $R_4.$
\end{corollary}

\begin{proof}
The modules of the set $\mathcal{M}$ depends on $b,c,d$ and
$\varepsilon,$ thus there are $3\cdot 3\cdot 3\cdot 2\cdot 2=108$
elements in this set. Since these modules are isomorphic in
couples, we obtain $54$ isomorphism classes which have the
representatives in the set $\mathcal{M}.$ The modules
$\mbox{Coker\ }\eta(a,b,c,\epsilon)$ depends on $\epsilon$ and on
the permutation $(a\ b\ c)$ of the cubic roots of $-1$ and the
modules $\mbox{Coker\ }\vartheta(a,b,c)$  are determined by the
permutation $(a\ b\ c)$ of the roots of $-1$, thus we get $6\cdot
2+6=18$ isomorphism classes which have the representatives in the
set $\mathcal{N}.$
\end{proof}

\thanks
We owe thanks to V. Vuletescu who used some results from algebraic
geometry to tell us in advance the number of non-isomorphic Ulrich
modules given by our Corollary \ref{victor}.

\end{document}